\documentstyle[12pt]{article}
\begin{document}

\begin{center}
\Large \bf Birationally rigid Fano hypersurfaces
\\ with isolated singularities
\end{center}
\vspace{1cm}

\centerline{\large  A.V.Pukhlikov}
\vspace{1cm}

\parshape=1
3cm 10cm
\noindent
{\small
\quad\quad\quad
\quad\quad\quad\quad \quad\quad\quad {\bf }\newline
It is proved that a general Fano hypersurface
$V=V_M\subset{\bf P}^M$ of index 1 with isolated singularities
of general position is birationally rigid. Therefore
it cannot be fibered into uniruled varieties of a smaller
dimension by a rational map and any ${\bf Q}$-Fano
variety $V'$ with Picard number 1 which is birational
to V is actually isomorphic to $V$. In particular, $V$ is
non-rational. The group of birational self-maps of $V$
is either $\{1\}$ or ${\bf Z}/2{\bf Z}$, depending on
whether $V$ has a terminal point of the maximal possible
multiplicity $M-2$. The proof is based upon the method 
of maximal singularities and the techniques of hypertangent systems
combined with the Shokurov connectedness principle.}
\vspace{1cm}

\centerline{}

\noindent
0. Introduction

0.1. Birationally rigid varieties

0.2. Regular hypersurfaces

0.3. The main result

0.4. Earlier results on singular Fano varieties

0.5. Acknowledgements

\noindent
1. Start of the proof

1.1. Maximal singularities

1.2. Isolated singular point

1.3. The crucial fact

\noindent
2. Infinitely near maximal singularities

2.1. Resolution of a maximal singularity

2.2. Simple examples

2.3. Subvarieties of codimension 2

2.4. The case $\mu\geq 5$

2.5. The harder cases $\mu=3$ and $4$

\noindent
3. The technique of counting multiplicities

3.1. A sequence of blow ups

3.2. The self-intersection of a linear system

3.3. Proof of Proposition 6

\section*{Introduction}
\subsection{Birationally rigid varieties}

In this paper we work over the field ${\bf C}$ of complex
numbers. Recall that a Fano variety $X$ of dimension
$\geq 3$ with ${\bf Q}$-factorial terminal singularities,
$\mathop{\rm rk}\mathop{\rm Pic} X=1$ is said to be
{\it birationally superrigid}, if for each birational map
$$
\chi\colon X-\,-\,\to X'
$$ 
onto a variety $X'$ of the same
dimension, smooth in codimension one, and each linear system
$\Sigma'$ on $X'$, free in codimension 1 (that is,
$\mathop{\rm codim}\nolimits \mathop{\rm Bs}\Sigma'\geq 2$),
the inequality
\begin{equation}
\label{001}
c(\Sigma,X)\leq c(\Sigma',X')
\end{equation}
holds, where $\Sigma=(\chi^{-1})_*\Sigma'$ is the proper inverse image
of $\Sigma'$ on $X$ with respect to $\chi$, and $c(\Sigma,X)=c(D,X)$
stands for the {\it threshold of canonical adjunction}
$$c(D,X)=\sup\{b/a|b,a\in{\bf Z}_+\setminus \{0\},
|aD+bK_X|\neq\emptyset\}$$
$D\in\Sigma$, and similarly for $\Sigma'$, $X'$. 
$X$ is said to be {\it birationally rigid}, if for each
$X'$, $\chi$, $\Sigma'$ there exists a birational
self-map $\chi^*\in \mathop{\rm Bir}X$ such that the
triple $X'$, $\chi\circ\chi^*$, $\Sigma'$
satisfies the condition (\ref{001}).\par
The following fact is well-known. \par

{\bf Proposition 1.} {\it Assume that $X$ is rigid. Then:

{\rm (i)} $X$ can not be fibered into uniruled varieties by a
non-trivial rational map,

{\rm (ii)} if $\chi\colon X -\, -\, \to X'$ is a birational map
onto a Fano variety $X'$  with ${\bf Q}$-factorial terminal
singularities such that
$\mathop{\rm Pic} X'\otimes{\bf Q}={\bf Q}K_{X'}$,
then $X'$ is (biregularly) isomorphic to $X$. If $X$ is
superrigid, then $\chi$ itself is a (biregular) 
isomorphism. In particular, in the superrigid case the
groups of birational and biregular self-maps coincide:
$$
\mathop{\rm Bir}X=\mathop{\rm Aut}X.
$$

{\rm (iii)} $X$ is non-rational.}

\subsection{Regular hypersurfaces}
Let $W=W_m\subset{\bf P}^N$ be a hypersurface of degree $m\leq N$
in the $N$-dimensional complex projective space. For a point
$x\in W$ choose a system of affine coordinates $(z_1,\dots,z_N)$
on ${\bf C}^N\subset{\bf P}^N$ with the origin at $x$ and write
down the equation of the hypersurface $W$ as 
$$
f=q_1+q_2+\dots+q_m,
$$
where $q_i(z_*)$ are homogeneous polynomials of degree $i$.

{\bf Definition 1.} The hypersurface $W$ is {\it regular} at a
smooth point $x\in W$, if the sequence
$$
q_1,\dots,q_k,
$$
$k=\min\{m,N-1\}$ is regular in ${\cal O}_{x,{{\bf P}}^N}$, that
is, the system of equations 
$$
q_1=\dots=q_k=0
$$
defines in ${\bf P}^N$ an algebraic subset of codimension $k$.

A dimension count, similar to the arguments of [P3, Sec. 1],
shows that a general (in the sense of Zariski topology on
$H^0({\bf P}^n,{\cal O}_{{\bf P}^N}(m))$) hypersurface $W$ is
regular at each point.

Let 
$$
V=V_M\subset {\bf P}={\bf P}^M
$$ 
be a hypersurface of 
degree $M$, with at most isolated singularities,
$$
f=q_1+q_2+\dots+q_M
$$
its equation with respect to a system of affine coordinates
$(z_1,\dots,z_M)$ with the origin at $x\in V$. Let
$$
\mu=\min\{k\in{\bf Z}_+|q_k\not\equiv 0\}=\mathop{\rm mult}\nolimits_xV
$$
be the multiplicity of $V$ at the point $x$. Assume
that $M-2\geq\mu\geq 2$, that is, $x\in\mathop{\rm Sing } V$.

{\bf Definition 2.} The hypersurface $V$ is regular at the
point $x$, if the following conditions are satisfied:

(i) the sequence $q_{\mu},\dots,q_M$ is regular in
${\cal O}_{x,{{\bf P}}}$;

(ii) the hypersurface $T_xV=\{q_{\mu}=0\}\subset
{\bf T}={\bf P}(T_x{\bf P})\cong{\bf P}^{M-1}$ is smooth and
regular at each point $y\in T_xV$;

(iii) for $\mu=3,4$ and $M\geq 7$ for any point $y\in T_xV$
none of the irreducible components of the closed algebraic 
set
\begin{equation}
\label{01}
\{q_{\mu}=q_{\mu+1}=\dots=q_6=0\}\cap
T_y(T_xV)\subset{\bf T}
\end{equation}
is contained in the quadric hypersurface
\begin{equation}
\label{02}
T_y(T_y(T_xV)\cap T_xV)\subset{\bf T};
\end{equation}
for $\mu=3$, $M=6$ it is sufficient that this condition 
holds with $q_5$ instead of $q_6$ in (\ref{01}).

The condition (iii) should be explained, the more so that
we somewhat abuse our notations: the symbol $T_y(T_xV)$
stands for the hyperplane in ${\bf T}$, which is tangent
to $T_xV$ at the point $y$. Since the hypersurface $T_xV$
is regular, the intersection
$$
T_y(T_xV)\cap T_xV
$$
is a hypersurface in the hyperplane $T_y(T_xV)$ with an
isolated singular point of multiplicity 2. The closed set
(\ref{01}) has dimension $\geq 1$, so that we require that
none of its component is contained in the quadric (\ref{02}).
This condition can be formulated in a different way: the
intersection of the cycle
$$
T_xV\cap T_y(T_y(T_xV)\cap T_xV)
$$
with the complete intersection
$$
\{q_{\mu+1}=\dots=q_6=0\}
$$
is of codimension precisely $9-\mu$ in ${\bf T}$.

{\bf Proposition 2.} {\it Let ${\cal V}_{\mu}(x)\subset
{\bf P}(H^0({\bf P},{\cal O}_{{\bf P}}(M)))$ be the space
of hypersurfaces of degree $M\geq 5$, which have a singularity
of multiplicity $\mu$, $2\leq\mu\leq M-2$ at the fixed
point $x\in{\bf P}$. The general (in the sense of Zariski
topology) hypersurface $V\in{\cal V}_{\mu}(x)$ is regular at
each of its points.}

Obviously, for a general $V\in{\cal V}_{\mu}(x)$ we have
$\mathop{\rm Sing } V=\{x\}$. Let us point out the following
question: for which $k$-uples of integers
$(\mu_1,\dots,\mu_k)\in\{2,\dots,M-2\}^k$ there exists a
hypersurface $V$, which is regular at each of its points and
has $k$ points $x_1,\dots,x_k$ of multiplicities 
$\mu_1,\dots,\mu_k$, respectively? One can show [P5] that for 
$\mu_i\equiv 2$ regular hypersurfaces exist for $k\leq M+1$,
however it seems that the precise limit value of $k$ is
considerably higher.

\subsection{The main result}

The main result of the present paper is the following

{\bf Theorem.} {\it Assume that the hypersurface $V$ is 
regular at each point.} (i)  {\it If for any point
$x\in V$ the estimate 
$\mathop{\rm mult}\nolimits_xV\leq M-3$ holds, then $V$
is a birationally superrigid variety.} (ii) {\it If
$x\in V$ is (the only) singular point of multiplicity
$M-2$, then the projection from this point,
$$
\pi\colon V-\, -\, \to{\bf P}^{M-1}
$$
is of degree 2 and there exists a birational involution
(the Galois involution) $\tau\in\mathop{\rm Bir} V$, which
permutes the points in the fibers of $\pi$. The variety
$V$ is birationally rigid and the exact sequence
$$
1\to\mathop{\rm Aut} V\to\mathop{\rm Bir} V\to
\langle \tau\rangle={\bf Z}/2{\bf Z}\to 1.
$$
holds. For a general $V$ obviously $\mathop{\rm Aut} V=\{1\}$,
so that}
$\mathop{\rm Bir} V=\langle \tau\rangle\cong{\bf Z}/2{\bf Z}$.

\subsection{Earlier results on singular Fano varieties}
The first example of a birationally rigid  singular Fano 
3-fold was made
by the quartic $V=V_4\subset{\bf P}^4$ with a unique double
point of general position $x\in V$ [P1]. As in the case of
arbitrary dimension, the projection
$$
\pi\colon V-\, -\, \to {\bf P}^3
$$
from the point $x$ is of degree 2 and determines the Galois
involution $\tau_x\in\mathop{\rm Bir}V$. However, in dimension
three the group of birational self-maps is much bigger. There
are exactly 24 lines through the point $x\in V$ on $V$ (in the
case of general position), $L_1,\dots,L_{24}\subset V$. Let
$L=L_i$ be one of them. The projection
$$
\pi_L\colon V-\, -\, \to {\bf P}^2
$$
from this line fibers $V$ into elliptic curves. More exactly,
for a general point $p\in{\bf P}^2$ the curve 
$C_p=\pi^{-1}_L(p)$ is a plane cubic, passing through the
point $x$. Taking $x$ to be the zero of the group law on
$C_p$, we get a birational involution:
$$
\begin{array}{rccc}
\tau_L\colon & V & -\, -\, \to & V, \\
\tau_L|_{C_p}\colon & z & \mapsto & -z.
\end{array}
$$
Set $\tau_0=\tau_x$, $\tau_i=\tau_{L_i}$. The following fact
is true [P1]: 

{\it The variety $V$ is birationally rigid. The involutions
$\tau_i$, $i=0,1,\dots,24$, generate in $\mathop{\rm Bir} V$
a subgroup $B(V)$ of finite index, which is their free
product. The following exact sequence holds:
}
$$
1\to B(V)=\mathop{*}\limits^{24}_{i=0}\langle \tau_i\rangle
\to \mathop{\rm Bir}V \to \mathop{\rm Aut} V \to 1
$$
Here the action of $\mathop{\rm Aut}V$ on $B(V)$ is defined
in the obvious way.

In [C] Corti essentially simplified the proof, using the
Shokurov connectedness theorem [K] for exclusion of an infinitely
near maximal singularity over the point $x$. Somewhat later
Cheltsov noted that, in its turn, this argument of Corti's
can be simplified, if one applies Shokurov connectedness
to the exceptional divisor $E\subset V_0\to V$ of the blow
up of the point $x$, $E\cong {\bf P}^1\times{\bf P}^1$.
Namely, if the point $x$ is not maximal itself, but
there is an infinitely near maximal singularity over it,
then there is a linear system on $E$, say
$\Sigma_E$ (possibly, with fixed components), of curves
of type $(m,m)$ such that the log pair
$(E,\frac{1}{m}\Sigma_E)$ is not log canonical. But this
fact leads to a contradiction. In fact, this has
already been proved in [P1], see the proof of the
``Graph lemma''. 

This way of arguing is used in the present paper when
we consider a singular point of the maximal multiplicity
$M-2$.

Furthermore, in [P2] a series of birationally superrigid 
singular Fano varieties of arbitrary dimension was produced:
double spaces of index 1 with a double point of general
position. In [P5] singular Fano hypersurfaces 
$V=V_M\subset{\bf P}^M$ with non-degenerate double points
were proved to be birationally rigid. Finally, Corti and 
Mella [CM] considered a larger class of quartic 3-folds
with isolated double points.

In the paper [CPR] 95 families
of weighted Fano 3-fold hypersurfaces 
$V_d\subset{\bf P}(a_0=1,a_1,a_2,a_3,a_4)$, $d=a_1+\dots+a_4$
were proved to be birationally rigid (honestly speaking,
one should say 94 families, since the family number one
in this list is exactly the family of smooth quartics
$V\subset{\bf P}^4$, which were proved to be superrigid in [IM]
30 years ago, which made the starting point of the whole
rigidity theory). The weighted Fano hypersurfaces have
terminal factor-singularities. The present paper deals
with hypersurface singularities only.

\subsection{Acknowledgements}
This work was carried out during my stay at the University
of Bayreuth as a Humboldt Research Fellow. I would like to
thank Alexander von Humboldt-Stiftung for the financial
support of my research. I am also very grateful to Prof.
Th.Peternell and his colleagues at the Mathematics Institute
of the University of Bayreuth for the excellent conditions
of work, stimulating creative atmosphere and general
hospitality.

\section{Start of the proof}

We prove the theorem by means of the method of maximal
singularities, see [IM,P3,P4].
\subsection{Maximal singularities}
Following the traditional scheme of arguments, let us consider
a linear system $\Sigma\subset |nH|$ on $V$, where
$H\in \mathop{\rm Pic} V$ is the class of a hyperplane section.

The linear system $\Sigma$ is assumed to be moving (that is,
it has no fixed components). 

{\bf Definition 3.} A geometric discrete valuation
$\nu\in{\cal N}(V)$ is said to be a {\it maximal singularity}
of the linear system $\Sigma$, if the {\it Noether-Fano
inequality} holds:
$$
\nu(\Sigma)>n\cdot \mathop{\rm discrepancy}(\nu).
$$

If $V$ is not superrigid, then there exists a moving linear
system with a maximal singularity.

Set $B=\mathop{\rm centre}(\nu)\subset V$ to be the centre of
the maximal singularity, an irreducible subvariety of $V$.

{\bf Proposition 3.} $B=x\in V$ {\it is a singular point of the
hypersurface $V$ of multiplicity } $\mu\geq 3$.

{\bf Proof.} Assume that $B\not\subset\mathop{\rm Sing } V$.
By [P3, Sec. 3] we can assert that
$\mathop{\rm codim}\nolimits_VB\geq 3$
(otherwise take a projective curve $C\subset B$, 
$C\cap\mathop{\rm Sing } V=\emptyset$. For this curve we
have the estimate $\mathop{\rm mult}\nolimits_C\Sigma>n$.
Since $C$ is contained in the smooth part of $V$, the
arguments of [P3] work and give a contradiction.) 
Consequently, for the cycle 
$$
Z=(D_1\bullet D_2)
$$ 
of scheme-theoretic intersection of general divisors
$D_1,D_2\in\Sigma$ we get the
estimate
$$
\mathop{\rm mult}\nolimits_BZ>4n^2.
$$
Now there exists a smooth point $x\in B$. By the
regularity condition the arguments of [P3] give the
opposite estimate:
$$
\mathop{\rm mult}\nolimits_BZ\leq
\frac{4}{M}\mathop{\rm deg} Z=4n^2.
$$
A contradiction. Therefore, $B=x\in V$ is a singular
point. It was shown in [P5] that the case
$\mu=\mathop{\rm mult}\nolimits_x V=2$ is impossible.
Q.E.D. for the proposition.

\subsection{Isolated singular point}
Let $x\in V$ be a regular singular point of multiplicity
$\mu\geq 3$, $\nu\in {\cal N}(V)$ a maximal singularity of the
system $\Sigma\subset |nH|$, $x=\mathop{\rm centre}(\nu)$,
$$
\begin{array}{cccc}
\displaystyle
\varphi_0\colon & V_0 & \to & V \\
\displaystyle
            & \bigcup & &  \\
\displaystyle
            & E & \to & x
\end{array}
$$
the blow up of the singular point, $E=E_0$ the exceptional divisor. 
Since
$\mathop{\rm Pic} V_0={\bf Z} H\oplus{\bf Z} E$, for the strict
transform of the linear system $\Sigma$ on $V_0$ we get
$$
\Sigma^0\subset |nH-\nu E|.
$$
Recall that $E\subset{\bf P}^{M-1}$ is a regular hypersurface
of degree $\mu\geq 3$.

{\bf Proposition 4.} (i) {\it For $\mu\leq M-3$ the divisor $E$
cannot make a maximal singularity.}

(ii) {\it For $\mu=M-2$ the birational involution 
$\tau\in\mathop{\rm Bir} V$ is defined by a linear system
$|(M-1)H-ME|$ on $V_0$. If the point $x$ is maximal for the
system $\Sigma$, that is, $\nu_0>n$, then $\nu_0/n\leq M/(M-1)$
and $\tau_*\Sigma$ is a moving linear system on $V$, for
which the point $x$ is not maximal.}

{\bf Proof.} (i) Assume the converse: the point $x$ is
maximal for the system $\Sigma$. Then $\nu_0>(M-\mu-1)n\geq 2n$.
Let $D_1,D_2\in\Sigma$ be two general divisors. For the effective 
cycle $Z=(D_1\bullet D_2)$ of codimension 2 on $V$ we have the
estimate
$$
\mathop{\rm mult}\nolimits_xZ
\geq(M-\mu-1)^2n^2\mu>Mn^2=\mathop{\rm deg} Z
$$
(since $(M-\mu-1)^2\mu>M$), which is impossible. The contradiction
proves (i).

(ii) Consider $\tau$ as an element of 
$\mathop{\rm Bir} V_0\cong \mathop{\rm Bir} V$. It is easy to see
that outside an invariant closed subset of codimension 2 the
involution $\tau$ is biregular on $V_0$ and its action on
$\mathop{\rm Pic} V_0$ is given by the formulas
$$
\begin{array}{l}
\tau^* H=(M-1)H-ME,\\
\tau^* E=\mu H-(\mu+1)E.
\end{array}
$$
If the point $x$ is maximal for $\Sigma$, that is, $\nu_0>n$,
then
$$
\tau^*|nH-\nu_0E|\subset |(n(M-1)-\nu_0\mu)H-
(nM-\nu_0(\mu+1))E|,
$$
where $\mu=M-2$.Obviously, 
$nM-\nu_0(\mu+1)\leq n(M-1)-\nu_0\mu$, so that the point $x$
is no more maximal for the system $\tau_*\Sigma$.

Set $T=T_xV\cap V$: this is an irreducible divisor on $V$.
Obviously, 
$$
T\sim (M-2)H-(M-1)E.
$$ 
If the system $\Sigma$ is
moving then for a general divisor $D\in\Sigma$ the cycle
$Z=(T\bullet D)$ is effective, so that we get: 
$$
\begin{array}{ccc}
\nu_0\mathop{\rm mult}\nolimits_xT & 
\leq\mathop{\rm mult}\nolimits_xZ\leq & \mathop{\rm deg} Z \\
\parallel & & \parallel \\
\nu_0(M-1)(M-2) &   & M(M-2)n,
\end{array}
$$
whence $\nu_0/n \leq M/(M-1)$, as we claimed it to be. 
Q.E.D. for the proposition.

\subsection{The crucial fact}
Now the theorem follows from the following crucial fact.

{\bf Proposition 5.} {\it If the point $x\in V$ is not maximal
for the linear system $\Sigma$, then there exists no maximal
singularity $\nu\in{\cal N}(V)$ such that}
$x=\mathop{\rm centre}(\nu)$.

Indeed, Proposition 5 means that if the point $x\in V$ is not
maximal for the linear system $\Sigma$, then this system has
no maximal singularities at all. And this is exactly birational
(super)rigidity.

{\bf Proof of Proposition 5 for $\mu=M-2$.} Since the point $x$
is not maximal, the existence of a maximal singularity over the
point $x$ implies the existence of a maximal singularity for
the system $\Sigma^0$: in terms of the log-minimal model
program, the pair $(V_0,\frac{1}{n}\Sigma^0)$ is not
canonical on $E$. By Shokurov connectedness theorem [K],
the pair $(E,\frac{1}{n}\Sigma^0|_E)$ is not log-canonical.
However, this is impossible [Ch]: the set  $Y\subset E$
where the pair $(E,\frac{1}{n}\Sigma^0|_E)$ is not
log-canonical, cannot be of positive dimension by [P3] and
cannot be purely zero-dimensional by [Ch](recall that the
linear system $\Sigma^0|_E$ is cut out on $E$ by hypersurfaces
of degree $\nu_0\leq n$, since we assumed that the point $x$
is not maximal, so that the pair
$$
(E,\frac{1}{\nu_0}\Sigma^0|_E)
$$
is also not log-canonical). This contradiction completes the
proof for $\mu=M-2$.

The arguments of [Ch] extend the arguments of [ChPk].

\section{Infinitely near maximal singularities}

\subsection{Resolution of a maximal singularity}
Recall the standard constructions and notations [P3,P4]. Let
$$
\begin{array}{cccc}
\displaystyle
\varphi_{i,i-1}: & V_i & \to & V_{i-1} \\
\displaystyle
            & \bigcup & & \bigcup \\
\displaystyle
            & E_i & \to & B_{i-1},
\end{array}
$$
$i=1,\dots,K$, be the resolution of a valuation 
$\nu\in{\cal N}(V)$, which is maximal for $\Sigma$. Here the
first $L$ blow ups correspond to the cycles $B_{i-1}$ of
codimension $\geq 3$ (the lower part), whereas the following
$K-L$ blow ups correspond to the cycles $B_{i-1}$ of codimension
2 (the upper part;  it is possible that $K=L$ and the upper
part is empty). Set $p_i=p_{Ki}$ to be the number of paths from
$E_K$ to $E_i$, $i=0,\dots,K$, in the oriented graph $\Gamma$
of the valuation $\nu$ (see [IM,P3,P4]). Set 
$\delta_i=\mathop{\rm codim}\nolimits B_{i-1}-1$,
$i=1,\dots,K$, $\delta_0=M-\mu-1$.

For an irreducible subvariety $Y\subset V$ of codimension 2
set 
$$
m(Y)=\mathop{\rm mult}\nolimits_xY, \quad
m_i(Y)=\mathop{\rm mult}\nolimits_{B_{i-1}}Y^{i-1},
$$
$i=1,\dots,L$, where the upper index $j$ means that we take
the strict transform of the subvariety on $V_j$. The following
statement makes the technical base of the proof.

{\bf Proposition 6.} {\it If $\nu\in{\cal N}(V)$ is a maximal
singularity of the system $\Sigma$, $\mathop{\rm centre}(\nu)=x$,
then there exists an irreducible subvariety $Y\subset V$ of
codimension 2, which satisfies the following estimate}:
\begin{equation}
\label{1}
\frac{2}{\mu}p_0m(Y)+
\sum^L_{i=1}p_im_i(Y)>
\frac{\left(\sum\limits^K_{i=0}p_i\delta_i\right)^2
}{\frac{1}{2}p_0+\sum\limits^K_{i=1}p_i}
\frac{\mathop{\rm deg} Y}{M}.
\end{equation}

{\bf Proof} is given below in Sec. 3.

{\bf Remark.} As it was shown in [P3,P4], it is possible
to ``correct'' the coefficients $p_i$ in such a way that
the estimate
\begin{equation}
\label{2}
p_0\leq\sum^L_{i=1}p_i
\end{equation}
holds, if only $L\geq 1$. In order to do this, it is sufficient
to erase in the graph $\Gamma$ the arrows connecting 
$E_i$, $i\geq L+1$, with $E=E_0$, if there are such arrows
(otherwise there is nothing to prove). After this operation the 
Noether-Fano inequality becomes stronger, whereas the proof
of Proposition 6 still holds. In what follows, if $L\geq 1$,
then we assume that (\ref{2}) is true without special
comments.

Fix an irreducible subvariety $Y\subset V$ of codimension 2,
which satisfies the estimate (\ref{1}). Our aim is to get a
contradiction and thus to show that the initial assumption
that there is a maximal singularity $\nu\in{\cal N}(V)$ with
the centre at the point $x$ is wrong. Birational rigidity of
$V$ would be an immediate implication of that.
 
\subsection{Simple examples}

{\bf Proposition 7.} $L\geq 1$.

{\bf Proof.} Assume that $L=0$. The estimate (\ref{1}) takes
the form of the inequality
$$
m(Y)>
\frac{\mu}{2}\frac{\left((M-\mu-1)p_0+\Sigma_u\right)^2
}{p_0(\frac{1}{2}p_0+\Sigma_u)}\frac{\mathop{\rm deg} Y}{M}.
$$
Here for convenience $\Sigma_u=\sum\limits^K_{i=L+1}p_i=
\sum\limits^K_{i=1}p_i$.
By the definition of the integers $p_i$ we get an obvious
estimate
$$
p_0\leq \Sigma_u.
$$
It is easy to check that for each $s,t$ the following
inequality holds:
$$
\frac{(2s+t)^2}{2s(\frac{s}{2}+t)}\geq 3,
$$
whence we get the estimate
$$
\mathop{\rm mult}\nolimits_xY>\frac{3\mu}{M}\mathop{\rm deg} Y.
$$

Let 
$$
f_i=q_{\mu}+\dots+q_i,
$$ 
$\mu\leq i\leq M$, denote the left
segment of the equation of the hypersurface $V$.

{\bf Definition 4.} The linear system
$$
\Lambda_i=|\sum^i_{j=\mu}f_js_{i-j}|_V,
$$
where $s_k(z_*)$ stands for an arbitrary homogeneous polynomial
of degree $k$ in $z_*$, is called the $i$-{\it th hypertangent
linear system} at the point $x$.

Obviously, for any divisor $D\in\Lambda_i$ we get
\begin{equation}
\label{3}
\frac{\mathop{\rm mult}\nolimits_x}{\mathop{\rm deg}} D
\geq \frac{i+1}{i}\frac{\mu}{M}.
\end{equation}
By the regularity condition
\begin{equation}
\label{4}
\mathop{\rm codim}\nolimits_V\mathop{\rm Bs} \Lambda_i=i-\mu+1,
\end{equation}
for $i=\mu,\dots,M-1$. Now let 
$D_{\mu},D_{\mu+1},\dots,D_{M-1}$ be general divisors of
the hypertangent linear systems 
$\Lambda_{\mu},\dots,\Lambda_{M-1}$, respectively. It is easy
to see that by (\ref{4}) the set-theoretic intersection
$$
Y\cap D_{\mu+2}\cap D_{\mu+3}\cap \dots\cap D_{M-1}
$$
is of pure codimension $M-\mu$ in $V$. Consider the
effective cycle
$$
Y^*=(Y\bullet D_{\mu+2}\bullet 
D_{\mu+3}\bullet \dots\bullet D_{M-1})
$$
of the corresponding scheme-theoretic intersection. By (\ref{3})
we get the estimate
$$
\frac{\mathop{\rm mult}\nolimits_x}{\mathop{\rm deg}} Y^*\geq
\frac{3\mu}{M}\cdot\frac{\mu+3}{\mu+2}
\cdot\cdots\cdot\frac{M}{M-1}=\frac{3\mu}{\mu+2}>1,
$$
which is impossible. This contradiction proves Proposition 7.

Set 
$$
R=\{q_{\mu}=q_{\mu+1}=0\}\cap V.
$$ 
This is an irreducible
cycle of codimension 2 (by the regularity condition). For a 
general $V$ its strict transform $\tilde R$ on $V_0$ is
non-singular in a neighborhood of the exceptional divisor.

{\bf Proposition 8.} $Y\neq R$.

{\bf Proof.} The regularity condition implies that
$$
\mathop{\rm mult}\nolimits_xR=
\frac{\mu+2}{M}\mathop{\rm deg} R, \quad
\mathop{\rm mult}\nolimits_{B_{i-1}}R^{i-1}\leq 1.
$$
Now if $Y=R$, then by (\ref{1}) we get
$$
2\frac{\mu+2}{\mu}p_0+\Sigma_l>
\frac{(2p_0+2\Sigma_l+\Sigma_u)^2}{\frac12p_0+
\Sigma_l+\Sigma_u},
$$
where we set for convenience $\Sigma_l=\sum\limits^L_{i=1}p_i$.
Elementary computations show that this inequality is false.
Q.E.D. for the proposition.

The two examples, considered above (Propositions 7 and 8) are
the simplest cases.

Let us study the general case.

\subsection{Subvarieties of codimension 2}
Let $y\in E\subset{\bf P}^{M-1}$ be an arbitrary point
on the exceptional divisor. The regularity condition
gives that $T_y^{(1)}=T_yE\cap E$ is a hypersurface
(of degree $\mu$) in the hyperplane 
$T_yE\cong{\bf P}^{M-2}$, with the point $y$ as an isolated
quadric singularity. Set
$T=T_y(T^{(1)}_y)\cap T^{(1)}_y$. This is an irreducible
cycle of codimension 2 on $E$. Obviously,
$\mathop{\rm deg} T=2\mu$, 
$\mathop{\rm mult}\nolimits_yT=6$.

{\bf Lemma 1.} {\it Let $W\neq T$ be an irreducible subvariety
of codimension 2 on $E$. The following estimate holds for}
$\mu\geq 4$:
$$
\frac{\mathop{\rm mult}\nolimits_y}{\mathop{\rm deg}} W
\leq\frac{8}{3\mu}.
$$

{\bf Proof.} Apply the technique of hypertangent systems to
$E\subset{\bf P}^{M-1}$. This is possible due to the 
regularity condition. More precisely, let
$(u_1,\dots,u_{M-1})$ be a system of linear coordinates on
${\bf P}^{M-1}$ with the origin at the point $y$,
$$
e(y)=\xi_1+\xi_2+\dots+\xi_{\mu}
$$
the equation of the hypersurface $E$, 
$e_i=\xi_1+\dots+\xi_i$ its left segment. Here $\xi_i$ are
homogeneous of degree $i$ in $u_*$. Set
$$
\Delta_i=|\sum^i_{j=1}e_js_{i-j}|_E,
$$
$i=1,\dots,\mu-1$, where $s_k$ is an arbitrary
homogeneous polynomial of degree $k$. Obviously, by the
regularity condition
$$
\mathop{\rm codim}\nolimits\mathop{\rm Bs}\Delta_i\geq i,
$$
so that for a general divisor $D_i\in\Delta_i$ and an
arbitrary subvariety $B\subset E$ of codimension $i-1$
we get $B\not\subset D_i$. Since 
$W\neq T=\mathop{\rm Bs}\Delta_2$, we get
$$
\mathop{\rm codim}\nolimits W\cap D_2=3,
$$
so that $(W\bullet D_2)$ is an effective cycle of
codimension 3 on $E$. Therefore
$$
W\cap D_2\cap D_4\cap D_5 \cap\dots\cap D_{\mu-1}
$$
is of codimension $\mu-1$ on $E$, so that
$$
W^*=(W\bullet D_2\bullet D_4\bullet D_5 \bullet\dots\bullet D_{\mu-1})
$$
is an effective cycle on $E$. We obtain the estimate
$$
\begin{array}{cc}
\displaystyle
1\geq\frac{\mathop{\rm mult}\nolimits_y}{\mathop{\rm deg}} W^*\geq 
\frac{\mathop{\rm mult}\nolimits_y}{\mathop{\rm deg}} W\cdot &
\displaystyle
\left(\frac32\cdot\frac54\cdot\,\cdots\,\cdot\frac{\mu}{\mu-1}\right)\\
\\
\displaystyle  & \displaystyle \parallel \\  \\
\displaystyle
  & 
\displaystyle
\frac{3\mu}{8},
\end{array}
$$
which immediately implies the lemma.

{\bf Lemma 2.} {\it Let $\mu=3$. For any irreducible
subvariety $W\neq T$ of codimension 2 we get}
$$
\frac{\mathop{\rm mult}\nolimits_y}{\mathop{\rm deg}} W
\leq\frac23.
$$

{\bf Proof.} In the notations of the proof of the 
previous lemma $\mathop{\rm codim}\nolimits(W\cap D_2)=3$,
so that $(W\bullet D_2)$ is an effective cycle of
codimension 3 and
$$
1\geq\frac{\mathop{\rm mult}\nolimits_y}{\mathop{\rm deg}} 
(W\bullet D_2)\geq 
\frac{\mathop{\rm mult}\nolimits_y}{\mathop{\rm deg}} W
\cdot\frac32,
$$
which is what we need. Q.E.D.

Let $Y^0\subset V_0$ be the strict transform of the subvariety
$Y$ and $(Y^0\bullet E)={\bf P}(T_xY)$ its projectivized
tangent cone at $x$. For the effective cycle
$(Y^0\bullet E)$ of codimension 2 on $E$ we get the
presentation
\begin{equation}
\label{5}
(Y^0\bullet E)=aT+W,
\end{equation}
where $a\in{\bf Z}_+$ and the effective cycle $W$ does not
contain $T$ as a component.

{\bf Lemma 3.} $a\geq 1$.

{\bf Proof.} Assume the converse: $a=0$. Consider first
the case $\mu\geq 4$. By Lemma 1 we get
$$
\frac{m_1(Y)}{m(Y)}\leq
\frac{\mathop{\rm mult}\nolimits_y(Y^0\bullet 
E)}{\mathop{\rm deg}(Y^0\bullet E)}\leq
\frac{8}{3\mu},
$$
where $y\in B_0$ is an arbitrary point (since
$m(Y)=\mathop{\rm mult}\nolimits_xY=
\mathop{\rm deg}(Y^0\bullet E)$ and 
$m_1(Y)\leq\mathop{\rm mult}\nolimits_yY^0
\leq\mathop{\rm mult}\nolimits_y(Y^0\bullet E)$).
Thus taking into account the inequality
$m_i(Y)\leq m_1(Y)$, we may replace $m_i(Y)$ in (\ref{1})
by $(8/3\mu)m(Y)$ for $i\geq 1$. Now from
(\ref{1}) we obtain
$$
\begin{array}{ccc}
\displaystyle
\mathop{\rm mult}\nolimits_xY> & \displaystyle
\frac{(2p_0+2\Sigma_l+\Sigma_u)^2}{
(\frac{2}{\mu}p_0+\frac{8}{3\mu}\Sigma_l)(\frac12 p_0+
\Sigma_l+\Sigma_u)} & \displaystyle \cdot\quad
\frac{\mathop{\rm deg} Y}{M} \\ \\ \displaystyle
   &  \displaystyle \parallel &   \\ \\
   & \displaystyle \mu\frac{
4p_0^2+4p_0(2\Sigma_l+\Sigma_u)+4\Sigma_l(\Sigma_l+\Sigma_u)
+\Sigma^2_u}{
p^2_0+2p_0(\frac53\Sigma_l+\Sigma_u)+\frac83\Sigma_l
(\Sigma_l+\Sigma_u)}, &
\end{array}
$$
whence we get finally that
\begin{equation}
\label{6}
\frac{\mathop{\rm mult}\nolimits_x}{\mathop{\rm deg}} Y>
\frac{3}{2M}\mu.
\end{equation}
On the other hand, arguing as in the proof of Lemma 1,
we see that for general divisors $D_i\in\Lambda_i$ the
set-theoretic intersection
$$
Y\cap D_{\mu+2}\cap D_{\mu+3}\cap \dots\cap D_{M-1}
$$
is of codimension precisely $M-\mu$, so that taking the
effective cycle
$$
Y^*=(Y\bullet D_{\mu+2}\bullet \dots\bullet D_{M-1}),
$$
we get the estimate
\begin{equation}
\label{7}
\frac{\mathop{\rm mult}\nolimits_x}{\mathop{\rm deg}} Y
\leq\frac{\mu+2}{M}.
\end{equation}
Comparing this inequality with (\ref{6}), we see that
$$
\frac32\mu<\mu+2,
$$
so that $\mu<4$: a contradiction.

Now assume that $\mu=3$. By Lemma 2, in this case
$$
\frac{m_1(Y)}{m(Y)}\leq\frac{2}{\mu},
$$
so that, arguing as above, we get the estimate
$$
\frac{\mathop{\rm mult}\nolimits_x}{\mathop{\rm deg}} Y> 3 
\frac{
4p_0^2+4p_0(2\Sigma_l+\Sigma_u)+4\Sigma_l(\Sigma_l+\Sigma_u)
+\Sigma^2_u}{
p^2_0+p_0(3\Sigma_l+2\Sigma_u)+2\Sigma_l
(\Sigma_l+\Sigma_u)}\cdot\frac{1}{M},
$$
so that
$$
\frac{\mathop{\rm mult}\nolimits_x}{\mathop{\rm deg}} Y>
\frac{6}{M}.
$$
The estimate (\ref{7}) is true for $\mu=3$, either, so that
we get a contradiction: $6<\mu+2=5$. Q.E.D. for Lemma 3.

{\bf Corollary 1(from Lemma 3).}  $Y\not\subset T_xV$.

{\bf Proof.}  Assume the converse: $Y\subset T_xV$. Then we
get $Y\subset T_xV\cap V$ and therefore 
$T_xY\subset T_x(T_xV\cap V)$. However, this is impossible,
since ${\bf P}(T_xY)$ contains the subvariety $T$ as a component,
whereas
$$
{\bf P}(T_x(T_xV\cap V))=\{q_{\mu}=q_{\mu+1}=0\}
\subset{\bf P}^{M-1}
$$
does not contain $T$ by the regularity condition. Q.E.D. for the
corollary.

\subsection{The case $\mu\geq 5$}
Now let us prove Proposition 5 for $\mu\geq 5$. For any 
irreducible subvariety $W\subset E$ of codimension 2 by Lemma 1
we get
$$
\frac{\mathop{\rm mult}\nolimits_y}{\mathop{\rm deg}} W\leq 
\frac{3}{\mu},
$$
where the equality is attained at $W=T$ only. Arguing as above,
we get from (\ref{1}):
\begin{equation}
\label{8}
\frac{\mathop{\rm mult}\nolimits_x}{\mathop{\rm deg}} Y>\mu \frac{
(2p_0+2\Sigma_l+\Sigma_u)^2}{M(2p_0+3\Sigma_l)(\frac12
p_0+\Sigma_l+\Sigma_u)}\geq
\frac{4}{3M}\mu.
\end{equation}
Since $Y\not\subset T_xV$, the intersection $Y\cap D_{\mu}$
is of codimension 3, so that by the regularity condition
$$
Y^*=(Y\bullet D_{\mu}\bullet D_{\mu+3}\bullet \dots\bullet D_{M-1})
$$
is an effective cycle of codimension $M-\mu$. Now we get
\begin{equation}
\label{9}
\begin{array}{ccc}
\displaystyle 
1\geq\frac{\mathop{\rm mult}\nolimits_x}{\mathop{\rm deg}} Y^*
\geq\frac{\mathop{\rm mult}\nolimits_x}{\mathop{\rm deg}} Y & \cdot &
\displaystyle
\left(
\frac{\mu+1}{\mu}\cdot\frac{\mu+4}{\mu+3}\cdot\cdots\cdot
\frac{M}{M-1}\right) \\   \\
  &   & \displaystyle \parallel \\ \\
  &   & \displaystyle \frac{(\mu+1)M}{\mu(\mu+3)},
\end{array}
\end{equation}
so that combining (\ref{8}) and (\ref{9}) we obtain the estimate
$$
\frac{\mu(\mu+3)}{\mu+1}>\frac43 \mu,
$$
whence $\mu<5$: a contradiction. Proposition 5 is proved for
$\mu\geq 5$.

\subsection{The harder cases  $\mu=3$ and $4$}
There are two hardest cases left: $\mu=3$ and $\mu=4$. We will
do the second case in full detail. Here one should employ more
delicate arguments than those above.

Recall that $T_xY$ contains $T$ as a non-trivial component and
thus $Y\not\subset D_{\mu}$, as above. The more so  
$Y\not\subset D_{\mu+1}$ for a general divisor 
$D_{\mu+1}\in \Lambda_{\mu+1}$. However by the regularity
condition one can say more: the intersection
$$
T\cap\{q_{\mu+1}=q_{\mu+2}=0\}
$$
is of codimension 2 in $T$. In particular, the linear system
$\Lambda^0_{\mu+1}|_T$ has no fixed components. Thus none of
the components of the closed algebraic set $D^0_{\mu+1}\cap T$
is contained in the support of the cycle $W$ (\ref{5}). 
Set $Y_{\mu+1}=(Y\bullet D_{\mu+1})$. This is an effective
cycle of codimension 3 on $V$. We get the following
presentation:
$$
Y_{\mu+1}=Y^{\sharp}_{\mu+1}+Y^+_{\mu+1},
$$
where an irreducible component $X$ of the cycle $Y_{\mu+1}$
comes into $Y^{\sharp}_{\mu+1}$ (and does not come into
$Y^+_{\mu+1}$) when and only when its strict transform
$X^0\subset V_0$ contains an irreducible component of the set
$(D^0_{\mu+1}\cap T)$. By what was said above,
$$
(\tilde Y^{\sharp}_{\mu+1}\bullet E)=
a^{\sharp}(T\bullet D^0_{\mu+1})+(\sharp),
$$
here $a^{\sharp}\geq a\geq 1$. For the cycle $Y^+_{\mu+1}$
we get the estimate
\begin{equation}
\label{10}
\frac{\mathop{\rm mult}\nolimits_x}{\mathop{\rm deg}} 
Y^+_{\mu+1}\leq \frac{\mu+3}{M},
\end{equation}
which is obtained in the usual way. However, one can say
much more about the cycle $Y^{\sharp}_{\mu+1}$: by construction
$$
Y^{\sharp}_{\mu+1}\not\subset D_{\mu}!
$$
Consequently, $(Y^{\sharp}_{\mu+1}\bullet D_{\mu})$ is an
effective cycle of codimension 4, so that we get
\begin{equation}
\label{11}
\frac{\mathop{\rm mult}\nolimits_x}{\mathop{\rm deg}} 
Y^{\sharp}_{\mu+1}\leq
\frac{\mu}{\mu+1}
\frac{\mathop{\rm mult}\nolimits_x}{\mathop{\rm deg}} 
(Y^{\sharp}_{\mu+1}\bullet D_{\mu})\leq
\frac{\mu(\mu+4)}{(\mu+1)M}.
\end{equation}
Now set 
$$
\begin{array}{l}
d^{\sharp}=\mathop{\rm deg} Y^{\sharp}_{\mu+1},\\
d^+=\mathop{\rm deg} Y^+_{\mu+1}, \\  \\
b^{\sharp}=a \mathop{\rm deg} T,\\
b^+=\mathop{\rm deg} W, \\  \\
\mathop{\rm deg} (\tilde Y^{\sharp}_{\mu+1}\bullet E)=
b^{\sharp}(\mu+2)+\delta^{\sharp},\\
\mathop{\rm deg} (\tilde Y^+_{\mu+1}\bullet E)=\delta^+.
\end{array}
$$
We get a system of inequalities,
$$
\begin{array}{l}
(b^{\sharp}+ b^+)(\mu+2)\leq (\mu+2)b^{\sharp}+
\delta^{\sharp}+\delta^+, \\  \\
(\mu+2)b^{\sharp}+ \delta^{\sharp}\leq
d^{\sharp}\frac{\mu(\mu+4)}{M(\mu+1)}, \\  \\
\delta^+\leq d^+\frac{\mu+3}{M},
\end{array}
$$
where 
$$
d^++d^{\sharp}=(\mu+1)\mathop{\rm deg} Y,
$$
and
$$
b^++b^{\sharp}=\mathop{\rm mult}\nolimits_xY=m_0.
$$
Note first of all that since the inequality (\ref{11}) is
stronger than (\ref{10}), we may assume that 
$\delta^{\sharp}=0$: otherwise replace $\delta^+$ by
$\delta^++\delta^{\sharp}$, $\delta^{\sharp}$ by $0$,
$d^{\sharp}$ by 
$$
d^{\sharp}-\delta^{\sharp}\frac{M(\mu+1)}{\mu(\mu+4)}
$$
and $d^+$ by
$$
d^++\delta^{\sharp}\frac{M(\mu+1)}{\mu(\mu+4)}.
$$
All the inequalities above are still true since
$$
\delta^{\sharp}\leq\delta^{\sharp}
\frac{(\mu+1)(\mu+3)}{\mu(\mu+4)}.
$$
Furthermore, by Lemma 1 for $m_1(Y)$ we have the following
estimate:
$$
m_1(Y)\leq\frac{3}{\mu}b^{\sharp}+\frac{8}{3\mu}b^+.
$$
Taking into account (\ref{1}), this implies
$$
b^{\sharp}\left(
\frac{2}{\mu}p_0+
\frac{3}{\mu}\Sigma_l\right)+
b^+\left(
\frac{2}{\mu}p_0+
\frac{8}{3\mu}\Sigma_l\right)>
\frac{(2p_0+2\Sigma_l+\Sigma_u)^2}{(\frac12
p_0+\Sigma_l+\Sigma_u)} \frac{\mathop{\rm deg} Y}{M}.
$$
Using the estimates, obtained above, we get now
$$
\begin{array}{c}
\displaystyle
\frac{\mu(\mu+4)}{(\mu+1)(\mu+2)}d^{\sharp}
\left(
\frac{2}{\mu}p_0+
\frac{3}{\mu}\Sigma_l\right)+
\frac{\mu+3}{\mu+2}
d^+\left(
\frac{2}{\mu}p_0+
\frac{8}{3\mu}\Sigma_l\right)> \\ \\ \displaystyle
> 
\frac{(2p_0+2\Sigma_l+\Sigma_u)^2}{(\frac12
p_0+\Sigma_l+\Sigma_u)} \frac{(d^{\sharp}+d^+)}{\mu+1}.
\end{array}
$$
By linearity, either
$$
\frac{\mu+4}{\mu+2}>
\frac{(2p_0+2\Sigma_l+\Sigma_u)^2}{(2p_0+3\Sigma_l)
(\frac12
p_0+\Sigma_l+\Sigma_u)}\geq \frac43,
$$
whence $\mu<4$ --- a contradiction, or
$$
\frac{(\mu+3)(\mu+1)}{\mu(\mu+2)}>
\frac{(2p_0+2\Sigma_l+\Sigma_u)^2}{(2p_0+\frac83 \Sigma_l)
(\frac12
p_0+\Sigma_l+\Sigma_u)}\geq \frac32,
$$
whence $\mu^2<2\mu+6$ --- a contradiction again. The case
$\mu=4$ is completed.

Note that the estimates which we obtained above are
sufficient to exclude the case of a point of multiplicity
$\mu=3$ on the sextic 5-fold. If $M\geq 7$ and $\mu=3$, then
to prove Proposition 5, one should start with the cycle
$(Y\bullet D_{\mu+1})$, then look at those components of
this cycle which contain components of the cycle
$(T\cap D^0_{\mu+1})$. Then one should intersect these
components with $D_{\mu+2}$ and take those components of
the intersection which contain components of the cycle
$(T\cap D^0_{\mu+1}\cap D^0_{\mu+2})$. Finally, one should
intersect them with $D_{\mu}$ (this is still possible by
the regularity condition). The estimates, obtained by means
of these manipulations, are already strong enough to exclude
the case $\mu=3$. The corresponding computations are
rather tiresome and for this reason we do not give them
here in detail.

Q.E.D. for Proposition 5 and for the main theorem.

\section{The technique of counting multiplicities}
In this section we generalize the technique of counting
multiplicities [P3,P4] to certain classes of singularities.

\subsection{A sequence of blow ups}
Let $x\in X$ be a germ of an isolated terminal ${\bf Q}$-factorial
singularity and
$$
\begin{array}{cccc}
\displaystyle
\varphi_{i,i-1}: & X_i & \to & X_{i-1} \\
\displaystyle
            & \bigcup & & \bigcup \\
\displaystyle
            & E_i & \to & B_{i-1}
\end{array}
$$
$i=1,\dots,K$, a sequence of blow ups with centres
$B_{i-1}\subset X_{i-1}$, where $B_0=x$. Let 
$E_i=\varphi^{-1}_{i,i-1}(B_{i-1})\subset X_i$ be the exceptional
divisors. We assume that the following conditions hold:

(i) $\varphi_{i,i-1}(B_i)=B_{i-1}$, that is, $B_i\subset E_i$;

(ii) the exceptional divisors $E_i\subset X_i$ are irreducible,
reduced and $X_i$ is ${\bf Q}$-factorial over a general point
of the cycle $B_{i-1}$.

Set $\delta_i=\mathop{\rm codim} B_{i-1}-1$. Obviously, we get
$$
(F_i\cdot (-E_i)^{\delta_i})=\mu_i\geq 1,
$$
where 
$$
\mu_i=\mathop{\rm mult}\nolimits_{B_{i-1}}V_{i-1},
$$
$F_i$ is a fiber of the morphism 
$\varphi_{i,i-1}\colon E_i\to B_{i-1}$.

For a cycle $Y_i\subset X_i$ we denote by the symbol
$Y^j\subset X_j$ its strict transform when it is well
defined. On the set of exceptional divisors $\{E_i\}$ we
define a structure of an oriented graph in the usual way:
$$
E_j\to E_i \quad\mbox{or}\quad i\to j,
$$
if $j>i$ and $B_{j-1}\subset E^{j-1}_i$ [IM,P1-P5].
For $j\to i$ we set
$$
\beta_{i,j}=
\sup\limits_{Y\subset E_j}
\frac{\mathop{\rm mult}\nolimits_{B_{i-1}}Y^{i-1}}{\mathop{\rm 
deg} Y}\in {\bf R}_+,
$$
where $\sup$ is taken over all the prime divisors 
$Y\subset E_i$, covering $B_{j-1}$ (on the other hand,
if $\varphi_{j,j-1}(Y)\neq B_{j-1}$, then 
$\mathop{\rm mult}\nolimits_{B_{i-1}}Y^{i-1}=0$), and
$$
\mathop{\rm deg} Y=(Y\cdot F_j\cdot (-E_j)^{\delta_j-1})
$$
is the ``degree'' of the intersection $Y\cap F_t$, 
$F_t=\varphi^{-1}_{i,i-1}(s)$, $s\in B_{i-1}$ is a general
point. 

For a path $\pi\in P(i,j)$, connecting $i$ with $j$, we define
its {\it weight} to be 
$$
\beta(\pi)=\prod^k_{\alpha=1}\beta_{i_{\alpha},i_{\alpha-1}},
$$
where 
$\pi=\{i=i_k\to i_{k-1}\to\dots\to i_{\alpha}\to i_{\alpha-1}
\to\dots\to i_0=j\}$.
We define the coefficients $w_{i,j}$ by the formula
$$
w_{i,j}=\sum_{\pi\in P(i,j)}\beta(\pi), \quad w_{i,i}=1.
$$

{\bf Lemma 4.} {\it The following equality holds}
$$
w_{i,j}=\sum_{k\to j} w_{i,k}\beta_{k,j}.
$$

{\bf Proof.} Take the disjoint union
$$
P(i,j)=
\coprod_{k\to j} P(i,k)\circ\{k\to j\},
$$
where $\circ\{k\to j\}$ means the extension of a path
from $i$ to $k$ to a path from $i$ to $j$ by adding the
arrow $k\to j$. Now the claim of the lemma is obvious by
the definition of the numbers $w_{i,j}$.

\subsection{The self-intersection of a linear system}
Now take linear system $\Sigma$ on $X$ without fixed 
components, and set $\Sigma^i$ to be its strict transform on $X_i$,
$D\in\Sigma$ its general divisor. We get
$$
D^i=\varphi^{*}_{i,i-1}(D^{i-1})-\nu_i E_i,
$$
so that
$$
D^K=\varphi^{*}_{K,0}(D)-
\sum^K_{i=1}\nu_i \varphi^{*}_{K,i}E_i.
$$
Let $D_1,D_2\in\Sigma$ be two general divisors. Define a sequence
of cycles of codimension two on $X_i$, setting
$$
\begin{array}{l}
\displaystyle
D_1\bullet D_2=Z_0,\\
\displaystyle
D^1_1\bullet D^2_2=Z^1_0+Z_1,\\
\displaystyle
\dots,\\
\displaystyle
D^i_1\bullet D^i_2=
(D^{i-1}_1\bullet D^{i-1}_2)^i+Z_i,\\
\displaystyle
\dots,
\end{array}
$$
where $Z_i\subset E_i$. From this presentation we get for
$i\leq L$, where $L$ is defined by the condition
$\mathop{\rm codim} B_{i-1}\geq 3$ for $i\leq L$:
$$
D^i_1\bullet D^i_2=
Z^i_0+Z^i_1+\dots+Z^i_{i-1}+Z_i.
$$
For any $j>i,j\leq L$ set
$$
m_{i,j}=\mathop{\rm mult}\nolimits_{B_{j-1}}(Z^{j-1}_i).
$$
Set also 
$$
d_i=\mathop{\rm deg} Z_i=(Z_i\cdot F_i\cdot (-E_i)^{\delta_i-1}).
$$
We get the following system of equalities:
$$
\begin{array}{l}
\displaystyle
\mu_1\nu^2_1+d_1=m_{0,1},\\
\displaystyle
\mu_2\nu^2_2+d_2=m_{0,2}+m_{1,2},\\
\displaystyle
\vdots\\
\displaystyle
\mu_i\nu^2_i+d_i=m_{0,i}+\dots+m_{i-1,i},\\
\displaystyle
\vdots\\
\displaystyle
\mu_L\nu^2_L+d_L=m_{0,L}+\dots+m_{L-1,L}.
\end{array}
$$
Now
$$
d_L\geq
\sum^K_{i=L+1}\mu_i\nu^2_i.
$$
Multiply the $i$-th equation by $w_{L,i}$ and put them all
together. In the right-hand part for each $i\geq 1$
we get the expression
\begin{equation}
\label{31}
\sum^L_{j=i+1}w_{L,j}m_{i,j}=
\sum^L_{j\to i}w_{L,j}m_{i,j}.
\end{equation}
However, by the definition of the numbers $\beta_{j,i}$
we have the estimate
$$
m_{i,j}\leq \beta_{j,i}d_i,
$$
so the (\ref{31}) can be bounded from above by the number
$$
d_i\sum_{j\to i}w_{L,j}\beta_{j,i}=d_i w_{L,i}.
$$
In the left-hand part for each $i\geq 1$ we see 
$d_i w_{L,i}$, so that, throwing away all the 
$m_{i,*}$, $i\geq 1$, from the right-hand part and all the
$d_i$, $i\geq 1$, from the left-hand part, we get finally:
\begin{equation}
\label{32}
\sum^L_{j=1}w_{L,j}m_{0,j}\geq
\sum^L_{j=1}w_{L,j}\mu_j\nu^2_j+
\sum^K_{i=L+1}\mu_i\nu^2_i.
\end{equation}

\subsection{Proof of Proposition 6}
Let us come back to the singular point $x\in V$, considered in
the present paper. We obviously get
$$
w_{i,j}=1 \quad\mbox{for}\quad i,j\geq 1,
$$
where, in accordance with the notations, which we use in this
paper, the sequence of blow ups $\varphi_{i,i-1}$ starts with
$\varphi_0$, and not with $\varphi_{1,0}$. For any divisor
$Y\subset E$ and a point $y\in Y$ we get the estimate
$$
\frac{\mathop{\rm mult}\nolimits_y}{\mathop{\rm deg}} Y\leq
\frac{2}{\mu},
$$
where the equality is attained at the divisor $T_yE\cap E$ only.
Indeed, by the regularity condition for general divisors
$R_i\in\Delta_i$ of the hypertangent linear systems on $E$
we get for $Y\neq R_1=T_yE\cap E$: the intersection
$$
Y\cap R_1\cap R_3\cap\dots\cap R_{\mu-1}
$$
is of codimension precisely $\mu-1$ on $E$, whence it follows
that the cycle
$$
Y^*=(Y\bullet R_1\bullet R_3\bullet\dots\bullet
R_{\mu-1})
$$
is effective, so that
$$
\begin{array}{ccc}
\displaystyle
1\geq\frac{\mathop{\rm mult}\nolimits_y}{\mathop{\rm deg}} Y^*
\geq\frac{\mathop{\rm mult}\nolimits_y}{\mathop{\rm deg}} Y & \cdot &
\displaystyle
\left(
\frac{2}{1}\cdot\frac{4}{3}\cdot\,\cdots\,\cdot
\frac{\mu}{\nu-1}\right) \\ \\
  &   & \displaystyle \parallel \\ \\
  &   & \displaystyle \frac{2\mu}{3},
\end{array}
$$
and thus
$$
\frac{\mathop{\rm mult}\nolimits_y}{\mathop{\rm deg}} Y
\leq \frac{3}{2\mu}<\frac{2}{\mu}.
$$
Consequently, $\beta_{i,0}\leq 2/\mu$ for all $i\to 0$. Now we get
from (\ref{32}):
\begin{equation}
\label{33}
\frac{2}{\mu}p_{L,0}m_0+
\sum^L_{i=1}p_{L,i}m_{i}\geq
2p_0\nu_0^2+
\sum^L_{i=1}p_{L,i}\nu^2_i+
\sum^K_{i=L+1}\nu^2_i.
\end{equation}
By the definition of the integers $p_{i,j}$ the estimate
(\ref{33}) implies the inequality
\begin{equation}
\label{34}
\frac{2}{\mu}p_0m_0+
\sum^L_{i=1}p_im_{i}\geq
2p_0\nu_0^2+
\sum^K_{i=1}p_i\nu^2_i,
\end{equation}
where $p_i=p_{K,i}$. Let us get a lower bound for the
right-hand part of (\ref{34}). By the Noether-Fano
inequality we get
\begin{equation}
\label{35}
\sum^K_{i=0}p_i\nu_i>
n\left(\sum^K_{i=0}p_i\delta_i
\right).
\end{equation}
Since
$$
\inf\limits_{\sum^K_{i=0}p_i\nu_i=C}\{
2p_0\nu_0^2+\sum^K_{i=1}p_i\nu^2_i\}=
\frac{C^2}{\displaystyle \frac12 p_0+\sum^K_{i=1}p_i},
$$
we get finally (taking into consideration that
$\mathop{\rm deg} Z=Mn^2$, $Z=(D_1\bullet D_2)$):
$$
\frac{2}{\mu}p_0\mathop{\rm mult}\nolimits_x Z+
\sum^L_{i=1}p_i\mathop{\rm mult}\nolimits_{B_{i-1}}Z^{i-1}>
\frac{\left(\displaystyle \sum^K_{i=0}p_i\delta_i\right)^2
}{\displaystyle
\frac12 p_0+\sum^K_{i=1}p_i}\cdot\frac{\mathop{\rm deg} Z}{M}.
$$
It remains to note that this inequality is linear in $Z_*$.
Therefore, there exists an irreducible component $Y$ of this 
cycle which satisfies (\ref{1}).

Q.E.D. for Proposition 6.

\section*{References}

\noindent
[IM]  Iskovskikh V.A. and Manin Yu.I., Three-dimensional quartics and
counterexamples to the L{\" u}roth problem. --
Math. USSR Sb. {\bf 15.1}, 1971, 141-166.
\vspace{0.5cm}

\noindent
[C] Corti A., Singularities of linear systems and 3-fold
birational geometry.  In: ``Explicit Birational Geometry 
of Threefolds'' CUP 2000, 259-312.
\vspace{0.5cm}

\noindent
[CM] Corti A., Mella M. Birational geometry of terminal
quartic 3-folds. I, preprint, math.AG/0102096
\vspace{0.5cm}

\noindent
[CPR] Corti A., Pukhlikov A. and Reid M., Fano 3-fold
hypersurfaces. In: ``Explicit Birational Geometry of 
Threefolds'', CUP 2000, 175-258.
\vspace{0.5cm}

\noindent
[Ch] Cheltsov I.A., Log-canonical thresholds on hypersurfaces,
to appear in: Sbornik: Mathematics (2001).
\vspace{0.5cm}

\noindent
[ChPk] Cheltsov I.A. and Park J., Log-canonical thresholds 
and generalized Eckard points, preprint, math.AG/0003121.
\vspace{0.5cm}

\noindent
[K] Koll{\'a}r J., et al., Flips and Abundance for
Algebraic Threefolds, Asterisque 211, 1993.
\vspace{0.5cm}

\noindent
[P1] Pukhlikov A.V., Birational automorphisms of a
three-dimensional quartic with an elementary
singularity. Math. USSR Sbornik. {\bf 63} (1989), 457-482.
\vspace{0.5cm}

\noindent
[P2] Pukhlikov A.V., Birational automorphisms of double
spaces with singularities. J. Math. Sci. {\bf 85} (1997),
no. 4, 2128-2141.
\vspace{0.5cm}

\noindent
[P3] Pukhlikov A.V., Birational automorphisms of Fano
hypersurfaces, Invent. Math. {\bf 134} (1998), no. 2,
401-426.
\vspace{0.5cm}

\noindent
[P4] Pukhlikov A.V., Essentials of the method of maximal
singularities. In: ``Explicit Birational Geometry of 
Threefolds'', CUP 2000, 73-100.
\vspace{0.5cm}

\noindent
[P5] Pukhlikov A.V., Birationally rigid singular Fano
hypersurfaces. To appear in: J. Math. Sci. (2001).
\vspace{0.5cm}

\begin{flushleft}
Steklov Mathematics Institute\\
Gubkina 8\\
117966 Moscow\\
RUSSIA\\
e-mail: {\it pukh@mi.ras.ru}
\end{flushleft}

\end{document}